\documentclass[a4paper,11pt]{amsart}
\usepackage[english]{certus}
\usepackage{verbatim}

\DeclareUnicodeCharacter{2217}{*}

\title[Representation theory of topological full groups and paradoxicality]{Representation theory of topological full groups of \'etale groupoids and paradoxicality}
\author{Vadim Alekseev}
\address{Vadim Alekseev, Technische Universit\"{a}t Dresden, Fakult\"{a}t Mathematik, Institut f\"{u}r Geometrie, 01062 Dresden, Germany}
\email{vadim.alekseev@tu-dresden.de}

\author{Martin Finn--Sell}
\address{Martin Finn--Sell, Universit{\"a}t Wien, Fakult\"{a}t f\"{u}r Mathematik, Oskar-Morgenstern-Platz 1,   1090 Wien, \"{O}sterreich }
\email{martin.finn-sell@univie.ac.at}


\begin{document}
\onehalfspace

\maketitle

\begin{abstract}
We provide a unified treatment of several results concerning full groups of ample groupoids and paradoxical decompositions attached to them. This includes a criterion for the full group of an ample groupoid being amenable as well as comparison of its orbit, Koopman and groupoid-left-regular representations. Besides that, we unify several recent results about paradoxicality in semigroups and groupoids, relating embeddings of Thompson's group $V$ into full groups of ample étale groupoids.
\end{abstract}

\section{Introduction}

The worlds of semigroups and étale groupoids are known to be connected since several decades. Indeed, to every ample groupoid one can naturally attach its inverse semigroup of ample bisections; on the other hand, the work of Exel \cite{MR2419901} constructs several topological groupoids naturally attached to any inverse semigroup. Both inverse semigroups and étale groupoids have rich representation theory, thus connecting each of them to the realm of \Cst-algebras. 

Another object of fundamental importance which can be constructed from an ample étale groupoid is its full group; this provides a rich source of discrete groups with interesting properties. One of such properties which makes good sense for groups, inverse semigroups, groupoids and \Cst-algebras is amenability, and therefore systematic understanding connections between amenability of these objects is of significant importance. Since the discovery of amenability by von Neumann \cite{neumannZurAllgemeinenTheorie1929}, it is known that amenability is the opposite of \emph{paradoxicality} as defined by Banach and Tarski \cite{banachDecompositionEnsemblesPoints1924}.

This paper is devoted to a systematic treatment of amenability and representation theory of full groups of ample étale groupoids. Conceptualising the results of \cite{Grig-Dudko} and \cite{MR3627119} where the corresponding theory was developed for weakly branch groups, we prove
\begin{ThmA}[see Theorem \ref{Thm:irred}]
Let $\mc G$ be an ample effective Hausdorff groupoid, $F(\mc G)$ its topological full group, $\lambda$ the left regular representation of $\mc G$, $\kappa$ and $\rho_x$ the Koopman resp. orbit representation of $F(\mc G)$. Then for the representations $\pi_{\lambda}, \pi_{\kappa}$ and $\rho_{x}$ we have:
\begin{enumerate}
\item If $\mc G$ is minimal, then: \begin{enumerate} \item $\pi_{\lambda} \sim \rho_{x}$ for every $x \in X$;
\item $\rho_{x}$ is irreducible for each $x \in X$;
\item For $x,y \in X$, with $x$ and $y$ in distinct $\mc G$-orbits, $\rho_{x}$ and $\rho_{y}$ are not unitarily equivalent. 
\end{enumerate}
\item If $\mc G$ is topologically amenable, then $\pi_{\kappa} \prec \pi_{\lambda}$.
\item If $\mc G$ is minimal and topologically amenable, then $\pi_{\kappa} \sim \pi_{\lambda}$.
\end{enumerate}
\end{ThmA}

As a consequence, we obtain the following criterion for the full group of such groupoid to be amenable:
\begin{ThmA}[see Theorem \ref{Thm:amenable} and Theorem \ref{Thm:not-c-star-simple}]
Let $\mc G$ be a minimal ample Hausdorff effective groupoid. Then $F(\mc G)$ is amenable if and only if the following three conditions are satisfied:
\begin{enumerate}
\item $\mc G^{a}$ admits an invariant mean,
\item $\mc G$ is topologically amenable and
\item the rigid stabiliser $F(\mc G)_{(U)}$ is amenable for each clopen subset $U \subset \mc G^{(0)}$.
\end{enumerate}
Moreover, if the rigid stabiliser $F(\mc G)_{(U)}$ is amenable for every clopen subset $U$, then the group $F(\mc G)$ is not $C^{*}$-simple.
\end{ThmA}

Furthermore, combining the results of \cite{MR3448402},  \cite{MR3548134} and \cite{Bonicke-Li} we obtain a unified Tarski alternative for ample groupoids which in particular gives the following result:

\begin{ThmA}[Theorem \ref{Thm:generalgroupoid}]
Let $\mc G$ be a topologically amenable, ample, principal Hausdorff groupoid with compact unit space $X$ such that $\mc G^{a}$ satisfies $\mc D=\mc J$ and has an almost unperforated type semigroup. Then the following are equivalent:
\begin{enumerate}
\item There is a Tarski matrix of degree 2 over $\mc G^{a}$.
\item There is a clopen subset $E \subset X$ such that $E$ is $(2,1)$-paradoxical.
\item $F(\mc G)$ contains a subgroup isomorphic to Thompson's group $V$.
\end{enumerate}
\end{ThmA}

Applied to the coarse groupoid of a discrete metric space, this yields

\begin{ThmA}[Theorem \ref{Thm:wobbling}]
Let $X$ be a uniformly discrete metric space of bounded geometry with property A, then the following are equivalent:
\begin{enumerate}
\item $X$ is not supramenable,
\item Thompson's group $V$ embeds into the wobbling group $W(X)$.
\end{enumerate}
\end{ThmA}

\subsection*{Acknowledgements}
The authors would like to thank Rufus Willett for helpful discussions around the topics of the present paper. We also thank the anonymous referee whose comments and suggestions helped greatly to improve the text of the paper.

% \begin{comment}
% Paragraph about Thompson's group $V$, and the context of this small project.

% Paragraph about the Birget-Nekrasheyvch representation, which is a special case of a general representation attached to full groups of Cantor groupoids (or the tight representation of the group of units of a Boolean inverse monoid)

% The purpose of this note is to explore the broader context in which this representation sits, by considering more general circumstances of a representation of $V$ into a \Cst-algebra. The motivation for doing this is twofold, the first to explore how $V$ appears to be deeply connected with the paradoxicality of metric spaces (following the ideas of Lawson et al, also previously using work of ... and ...) as well as to identify potential \Cst-algebraic approaches to proving the Thompson's groups are \Cst-exact.

% A brief summary of our results, and the structure of the paper, are as follows.
% \end{comment}

\section{Inverse semigroups, groupoids and their \Cst-algebras}

\subsection{Inverse semigroups}
We consider the general text by Mark Lawson \cite{MR1694900} for a complete introduction to inverse semigroup theory, and recall the basic notions with references below. Recall that a \textit{semigroup} is a set $S$, together with an associative binary operation. If additionally it has a unit element, then we say it is a \textit{monoid}. 

\begin{Def}\label{Def:invsemi}
Let $S$ be a semigroup. We say $S$ is $inverse$ if there exists a unary operation $\ast:S \rightarrow S$ satisfying the following identities:
\begin{enumerate}
\item $(s^{*})^{*}=s$.
\item $s*$ is the unique element satisfyung $ss^{*}s=s$ and $s^{*}ss^{*}=s^{*}$ for all $s \in S$
\item $ef=fe$ for all idempotents $e,f \in S$.
\end{enumerate}
An element $e \in S$ satisfying $e^{2}=e$ is called an \textit{idempotent}, and the set of idempotents is denoted $E(S)$ -- in an inverse semigroup this is a (commutative) subsemigroup. Inverse semigroups also carry a natural partial order by setting $e\leq f$ iff $ef = e$. 
\end{Def}
\begin{Def}
Two elements $s$ and $t$ in an inverse semigroup $S$ (with zero) are \textit{compatible} if $s^*t$ and $st^*$ are elements of $E(S)$, and \textit{orthogonal} if $s^{*}t=st^{*}=0$. A finite set of elements $F \subset S$ is a compatible (resp. orthogonal) set if each pair of elements $s,t\in F$ are compatible (resp. orthogonal).
\end{Def}

A fundamental example of such an object is the \textit{symmetric inverse monoid} on any set $X$, denoted by $I(X)$. This is defined equipping the collection of all partial bijections of $X$ to itself equipped the composition defined on common intersections. The Wagner--Preston theorem (Section 1.5, Theorem 1 in \cite{MR1694900}), the analogue of Cayley's theorem for inverse semigroups, says that every abstract inverse semigroup can be realised as a subsemigroup of $I(X)$ for some set $X$. The natural order, $s \leqslant t$ then says that $t$ extends $s$ on a larger domain, and the compatibility condition says that $s$ and $t$ agree on the intersection of their domain, and that their inverses agree on the intersection of their ranges.

In this situation the function, $s\vee t$, defined by doing both $s$ and $t$ simultaneously on the union of their domains, is well defined and belongs to $I(X)$. A subsemigroup of $I(X)$ is \textit{distributive} if, given any compatible subset $F \subset S$, the join over $F$: $\bigvee_{s \in F}s$ also belongs to $S$.

\subsection{Groupoids}
For groupoids in the context in which we are studying them, useful references are \cite{MR2419901}, \cite{MR1799683}, \cite{MR1724106} and \cite{MR584266}. \textit{Topological groupoids} are elementary models for `noncommutative spaces' and appear throughout noncommutative geometry, index theory and operator algebras as a primary source of examples. 

\begin{Def}\label{def:grpoid2}
A \textit{groupoid} is a set $\mc G$ equipped with the following information:
\begin{enumerate}
\item A subset $\mc G^{(0)}$ consisting of the objects of $\mc G$, denote the inclusion map by $i: \mc G^{(0)}\hookrightarrow \mc G$. 
\item Two maps, $r$ and $s: \mc G  \rightarrow \mc G^{(0)}$ such that $r\circ i = s \circ i = Id$ .
\item An involution map $^{-1}:\mc G \rightarrow \mc G$ such that $s(g)=r(g^{-1})$.
\item A partial product $\mc G^{(2)} \rightarrow \mc G$ denoted $(g,h) \mapsto gh$, with $\mc G^{(2)}=\lbrace (g,h) \in \mc G \times \mc G | s(g)=r(h) \rbrace\subseteq \mc G\times \mc G$ being the set of pairs it is possible to compose.
\end{enumerate}
Moreover we ask the following:
\begin{itemize}
\item The product is associative where it is defined in the sense that for any pairs $(g,h),(h,k)\in \mc G^{(2)}$ we have $(gh)k $ and $g(hk)$ are defined and equal.
\item For all $g \in \mc G$ we have $r(g)g=gs(g)=g$.
\end{itemize}
\end{Def}

A groupoid is \textit{principal} if $(r,s): \mc G \rightarrow \mc G^{(0)} \times \mc G^{(0)}$ is injective. A groupoid $\mc G$ is a \textit{topological groupoid} if both $\mc G$ and $\mc G^{(0)}$ are topological spaces, and the maps $r,s, ^{-1}$ and the composition are all continuous. A Hausdorff, locally compact topological groupoid $\mc G$ is \textit{proper} if $(r,s)$ is a proper map and \textit{\'etale} or \textit{r-discrete} if the map $r$ is a local homeomorphism. When $\mc G$ is \'etale, $s$ and the product are also local homeomorphisms, and $\mc G^{(0)}$ is an open subset of $\mc G$ \cite[Section 3]{MR2419901}. A subset $Z\subset \mc G^{(0)}$ is called invariant if for all $g\in \mc G$ with $s(g) \in Z$ we have $r(g)\in Z$. If $x\in \mc G^{(0)}$, the orbit $\mc G x$ of $x$ is defined as the set $\{r(g)\mid g\in \mc G,\, s(g) = x\}$.

\begin{Def}
A \textit{bisection} in $\mc G$ is a subset $U \subset \mc G$ such that the range (or source) map $r|_U\colon U \to r(U)$ is a bijection. If $\mc G$ is \'etale, the set of open bisections $\mc G^{o}$ forms a basis for the topology. We say that an \'etale groupoid is \textit{ample} if the set of compact open bisections $\mc G^{a}$ is a basis for the topology of $\mc G$.% (this follows more along the lines of the definition given in Paterson \cite{MR1724106} or Renault \cite{MR584266}. 
\end{Def}

\paragraph{\textbf{Convention.}} From now on, all groupoids considered in this paper will be assumed ample and with compact base space $\mc G^{(0)}$.

\begin{Ex}
Let $G$ be a discrete group and let $X$ be a compact, Hausdorff topological $G$-space. Then the product $G\times X$ can be equipped with a groupoid structure that encodes the action of $G$ on $X$, as follows, with product given by $(g,x)(h,y)=(gh,y)$ whenever $x=hy$, inverse $(g,x)^{-1}=(g^{-1},gx)$ and source and range maps $s(g,x)=x$, $r(g,x)=gx$. We can topologise this by considering the sets $(g,U)= \lbrace (g,x) \mid x \in U \rbrace$, where $U$ is a open subset of $X$. This topological groupoid is denoted by $G\ltimes X$ and is called the \emph{action groupoid}. 
\end{Ex}

In the above example stabilisers may occur. To produce a principal groupoid, one uses the groupoid of germs construction.

\begin{Def}
Let $f,g$ be partial homeomorphisms of $X$ and let $x\in X$ be in the domain of both $f$ and $g$. Then $f$ and $g$ have the same germ at $x$, denoted $f \sim_{x} g$, if there is some neighbourhood $U$ of $x$ on which $f$ and $g$ agree.  
\end{Def}

Thus we can define, for any \'etale groupoid $\mc G$ a corresponding groupoid of germs $\Germs(\mc G)$ by considering the semigroup of compact open bisections $\mc G^a$ and letting $\Germs(\mc G)$ be the set of equivalence classes of germs of bisections. This inherits the product, inverse, range and source maps from the semigroup $\mc G^{a}$, and is a surjective groupoid image of $\mc G$, by mapping an element $\gamma \in \mc G$ to the germ of any bisection containing $\gamma$ at $s(\gamma)$. By the remark above, this map is surjective (as every element is contained in some bisection).% but not injective in general.\todo{add references}

$\Germs(\mc G)$ can be given a topology in the following way: for a clopen set $U \in \mc G^{(0)}=\Germs(\mc G)^{(0)}$, and an element in $A \in \mc G^{a}$ we can consider the sets $O_{A} = \lbrace [A_{x}] \mid x \in s(A) \rbrace$, and note that by declaring these sets be clopen when appropriate shows immediately that the map $\mc G \rightarrow \Germs(\mc G)$ defined above is open. 

\begin{Def}
An ample groupoid $\mc G$ is
\begin{itemize}
\item \textit{effective}, if for every non-identity $g \in \mc G$ and every bisection $A\in \mc G^{a}$ containing $g$, there is an element $h \in A$ such that $s(h)\not = r(h)$,
\item \textit{essentially principal} or \textit{topologically principal} if the set of points with trivial isotropy is dense in the unit space of $\mc G$,
\item \emph{minimal} if the only non-empty closed invariant subset of $\mc G^{(0)}$ is $\mc G^{(0)}$.
\end{itemize}
\end{Def}

The following is a standard result (see \cite[Proposition 2.1]{Nek-sgodo}, \cite{LL-Milan}, or \cite[Lemma 3.1]{MR3189105}):

\begin{Prop}\label{prop:nek}
If $\mc G$ is a topologically principal, Hausdorff ample groupoid, then it is effective. If $\mc G$ is second countable, Hausdorff effective groupoid, then $\mc G$ is topologically principal. \qed
\end{Prop}

\subsection{Boolean inverse monoids and full groups}

In this section we give specifics concerning the semigroup structure on bisections, and the reconstruction techniques that make noncommutative Stone duality work. 

The set of ample bisections $\mc G^{a}$ is a \textit{distibutive inverse semigroup} under composition of bisections -- if $A$ and $B$ are bisections, then the set product $AB$ is also a bisection, as is $A^{-1}$ and $A\cup B$ whenever $A$ and $B$ are compatible bisections. Note that the idempotent elements in $\mc G^{a}$ are clopen subsets of $\mc G^{(0)}$.

If $\mc G^{(0)}$ is compact, then $E(\mc G^{a})$ is a Boolean algebra, and $\mc G^{a}$ is a \textit{Boolean inverse monoid} \cite{MR3077869}, that is an distributive inverse monoid with a Boolean algebra of idempotents. These have been studied in the context of Stone dualities in the noncommutative setting \cite{MR3077869}

\begin{Def}
Let $\mc G$ be a second countable ample groupoid with compact base space. Then the group
\[
F(\mc G) = \{ U \in \mc G^a\mid r(U ) = s(U) = \mc G^{(0)}\}
\]
with composition inherited from $\mc G^a$ is called the topological full group of $\mc G$.
\end{Def}

The relationship of the full group to the original groupoid is described below.

\begin{Lemma}\label{Lem:germing}
Let $F(\mc G)$ be the topological full group of an ample topologically principal groupoid $\mc G$ with compact base space $\mc G^{(0)}$. Then $\Germs(F(\mc G)\curvearrowright \mc G^{(0)})$ is a closed subgroupoid of $\mc G$.
\end{Lemma}
\begin{proof}
The subgroupoid $\mc G^{(0)} \rtimes F(\mc G)$ is closed in $\mc G^{(0)} \rtimes \mc G^{a}$ (when $\mc G^{a}$ is given the discrete topology), and because the germ topology coincides with the quotient topology when $\mc G$ is Hausdorff, the quotient map is in fact closed.
\end{proof}

We can improve this precisely when the full group completely describes the local structure of $\mc G^{a}$ in the following sense. Indeed, by \cite[Theorem 2.22]{Lawson-minimal}, we can strengthen Lemma \ref{Lem:germing} to the following statement:

\begin{Prop}\label{Lem:germing2}
Let $F(\mc G)$ be the full group of an ample essentially principal minimal second countable groupoid $\mc G$ with compact base space $\mc G^{(0)}$. Then $\Germs(F(\mc G)\curvearrowright \mc G^{(0)})$ is homeomorphic to $\mc G$. \qed
\end{Prop}

\section{Representations of inverse semigroups, groupoids and their full groups}

In this section we recall the construction of groupoid \Cst-algebras. Following this, we illustrate that the groupoid representation associated to a full group (or certain subgroups of a full group) appears naturally as a reduction of the left regular representation of the associated Boolean inverse monoid of clopen bisections. This allows us to also formulate a clean connection between the rigid stabilisers of clopens and the groupoid representation being considered.
 
\subsection{Groupoids associated to inverse semigroups}

For any inverse semigroup $S$, it is possible to manufacture a variety of ample \'etale groupoids from $S$ that satisfy different purposes. We first recall some definitions, following the notation from \cite{MR2419901} and \cite{MR3548134}.

\begin{Def}
A \textit{filter} in $E(S)$ is a non-empty subset $\eta \subset E(S)$ such that:
\begin{enumerate}
\item $0 \not \in \eta$,
\item if $e,f \in \eta$, then $ef \in \eta$ and,
\item $e \in \eta$ and $e\leqslant f$, then $f \in \eta$.
\end{enumerate}
The set of filters is denoted $\widehat{E}(S)$, and can be viewed as a subspace of $\textbf{2}^{E(S)}$. For finite sets $X,Y \subset E(S)$, let
\begin{equation*}
U(X,Y) = \lbrace \eta \in \widehat{E}(S) \mid X \subset \eta, Y \cap \eta = \emptyset \rbrace.
\end{equation*}
The sets of this form are clopen and generate the topology on $\widehat{E}(S)$, as $X$ and $Y$ are varied over all finite subsets of $E(S)$. With this topology, the space $\widehat{E}(S)$  is called the0\emph{spectrum} of $E(S)$.

Recall that an \textit{ultrafilter} is a filter that is not properly contained in any other filter. The set of ultrafilters is denoted by $\widehat{E}_{\infty}(S)$, and as a subspace of $\widehat{E}(S)$ this \textit{may not be closed}. Let $\widehat{E}_{tight}(S)$ denote the closure of $\widehat{E}_{\infty}(S)$ in the topology of $\widehat{E}(S)$ -- when $E(S)$ is a Boolean algebra, we know immediately that $\widehat{E}_{tight}(S) = \widehat{E}_{\infty}(S)$ by Stone duality.
\end{Def}

The second definition we recall is that of an inverse semigroup action:

\begin{Def}
An \textit{action} of an inverse semigroup $S$ on a locally compact space $X$ is a semigroup homomorphism $\alpha: S \rightarrow I(X)$ such that:
\begin{enumerate}
\item $\alpha(s)$ is continuous for each $s \in S$,
\item the domain of $\alpha(s)$ is open for each $s \in S$, and
\item the union of the domains of the $\alpha(s),\,s\in S$ is equal to $X$.
\end{enumerate}
If $\alpha$ is an action of $S$ on $X$, then as with group actions, we will write $\alpha: S \curvearrowright X$. The above implies that $\alpha(s)^{-1}=\alpha(s^{*})$, and that each $\alpha(s)$ is a partial homemorphism of $X$. For each $e \in E(S)$, the map $\alpha(e)$ is the identity on some open subset $D_{e}^{\alpha}$, and one easily sees that the domain of $\alpha(s)$ is $D^{\alpha}_{s^{*}s}$, and the range is $D^{\alpha}_{ss^{*}}$, that is:
\begin{equation*}
\alpha(s): D^{\alpha}_{s^{*}s} \rightarrow D^{\alpha}_{ss^{*}}.
\end{equation*}
\end{Def}
We can now introduce the two groupoids we will associate to any inverse semigroup $S$, and recall some facts concerning them from the literature that we will need in the sequel.

\begin{Ex}
There is a natural action of $S$ on $\widehat{E}(S)$, and the \textit{universal groupoid} of $S$, denoted $G(S)$, is the transformation groupoid of this action. This was introduced by Paterson in \cite{MR1724106}, and appears in \cite{MR2419901} and has alternative descriptions, for instance found in % Lenz \cite{}\todo{add reference} and 
\cite{MR3077869}.
\end{Ex}

\begin{Ex}\label{ex:tight-reduction}
In addition to the universal groupoid, there is the \textit{tight reduction} $G_{tight}(S)$, which is the reduction of $G(S)$ to the invariant subspace $\widehat{E}_{tight}(S)$. This was introduced by Exel \cite{MR2419901}, and various of its properties have been studied in \cite{MR2534230} and \cite{MR2644910} -- an alternative description of this groupoid would be to consider the natural action of $S$ on $\widehat{E}_{tight}(S)$ and then take the transformation groupoid of it (following for instance the ideas of Exel \cite{MR2644910}).
\end{Ex}

In \cite{MR2644910}, Exel proved a reconstruction theorem, which is useful in this context, which states that if one begins with an ample groupoid $\mc G$, then $G_{tight}(\mc G^{a}) \cong \mc G$ via the procedure above applied to $S = \mc G^{a}$. More generally, this is an example of a \textit{noncommutative Stone duality} \cite{Lawson-minimal}, and has been studied in much more generality for classes of inverse semigroups by Lawson--Lenz \cite{MR3077869} -- it has even been improved to a categorical notion. For an inverse semigroup, we remark that the above procedure performed in a loop provides:
\begin{equation*}
G_{tight}(G_{tight}(S)^{a})\cong G_{tight}(S)
\end{equation*}
for any inverse semigroup $S$. 

\subsection{Inverse semigroup \Cst-algebras}

In this section we study from the perspective of the inverse monoid associated representations of the group of invertible elements, with the examples from the previous section in mind. First we recall the construction, and then discuss the notions of induction and restriction in this context.

Let $S$ be an inverse monoid. Following \cite[Section 2.1]{MR1724106}, we define the monoid algebra $\mb C S$ by extending the multiplication on $S$ linearly and the inversion antilinearly:
\[
 \left(\sum_{s\in S} {a_s s}\right)\left(\sum_{t\in S} {a_t t}\right) = \sum_{s,t\in S} a_s b_t\cdot st,
\]
\[
 \left(\sum_{s\in S} {a_s s}\right)^*  = \sum_{s\in S} \overline{a_s} s^*.
\]

This makes $\mb C S$ a $\ast$-algebra, and every element $s\in S$ is a partial isometry in $\mb C S$. Therefore, we can define the \emph{universal $C^*$-algebra of $S$} by completing $\mb C S$ with respect to the norm
 \[
  \norm{a} = \sup \norm{\pi(a)},
 \]
 where the supremum is taken over all $\ast$-representations of $\mb C S$ (the construction of the left regular representation in the next paragraph guarantees that the set of $\ast$-representations is non-empty). Similar to the case of discrete groups, $\ast$-representation of $S$ are in one-to-one correspondence with $\ast$-representations of $C^*(S):=\overline{\mathbb{C}S}^{\Vert . \Vert}$.

We denote by $\ell^2 S$ the Hilbert space with basis $S$ and let
\[
 \lambda\colon S\to \mb B(\ell^2 S)
\]
\[
 \lambda(s)\delta_t:= \begin{cases}
                     \delta_{st}, & tt^*\leqslant s^*s,\\
                     0, & \text{otherwise}
                    \end{cases}
\]
be the left regular representation of $S$; it extends to an injective representation of $\mb C S$ (in fact, by \cite{MR663865}, it extends to a injective representation of $\ell^{1}S$). Observe that if $G\subset S$ is a subgroup of $S$ that shares the unit of $S$, then $\lambda(G)$ consists of unitaries.
The \emph{reduced $C^*$-algebra of $S$} is then defined as $C^*_r(S):= \overline{\lambda(\mb C S)}\subset \mb B(\ell^2 S)$.

If $S$ has a zero element $o\in S$, then $I_o := \mb C \cdot o$ is an ideal in $\mb C S$ and we have
\[
 \lambda(o)\delta_s = 0,\quad s\neq o,\quad \lambda(s)\delta_o = \delta_o,\quad s\in S.
\]
We remark that $o$ gives rise to a central projection in any $\ast$-representation of $S$ and so if we consider the algebras:
\[
 C^*_o(S) = C^*(S) /I_o,
\]
\[
 C^*_{r,o}(S) = C^*_{r}(S)/ \lambda(I_o)
\]
we obtain a decomposition $C^{*}_{-}(S) \cong C^*_{-,o}(S)\oplus I_{o}$, where $-$ is a place holder for maximal or reduced respectively. We call $C^*_o(S)$ and $C^*_{r,o}(S)$ the \emph{truncated} universal and reduced $C^*$-algebras of $S$. 

Notice that in discussing properties like nuclearity or exactness it does not matter whether to work with the truncated or full version of inverse monoid $C^*$-algebras. It is also easy to see that $\ast$-representations of $C^*_{o}(S)$ correspond to $\ast$-representations of $S$ where $o$ acts as the zero operator.

\subsection{Groupoid \Cst-algebras and some representation theory}
Let $\mc G$ be a locally compact \'etale Hausdorff groupoid. Then the algebra of continuous, compactly supported functions
\begin{equation*}
C_{c}(\mc G) = \lbrace f: \mc G \rightarrow \mathbb{C} \mid \supp f \mbox{ compact} \rbrace
\end{equation*}
is a $*$-algebra under convolution and pointwise conjugation. It admits a maximal norm \cite{MR1724106,MR584266} satisfying the usual universal property. Denote the completion in this maximal norm by $C^{*}\mc G$ when there is no ambiguity. In addition to this maximal norm, there also a \textit{reduced norm}, defined as follows:
\begin{equation*}
\Vert f \Vert_{r}:= \sup_{x \in \mc G^{(0)}} \Vert \lambda_{x}(f) \Vert,
\end{equation*}
where $\lambda_{x}: C_{c}(\mc G) \rightarrow \mathfrak{B}(\ell^{2}(\mc Gx))$ is the natural left regular representation of $\mc G$ on the orbit of $x$, explicitly given by
\[
\lambda_x(f)(\delta_y) = \sum_{\substack{g\in \mc G\\ s(g) = y}} f(g)\delta_{r(g)},\quad y\in \mc G x.
\]

Denote the completion of $C_{c}(\mc G)$ in this norm by $C^{*}_{r}(\mc G)$, the reduced \Cst-algebra of $\mc G$. Note that by the universal property of $C^{*}(\mc G)$, there is a surjective quotient homomorphism $\lambda\colon C^{*}\mc G \rightarrow C^{*}_{r}(\mc G)$ called the left regular representation.

Given a closed subset $C \subset \mc G^{(0)}$, we define the restriction groupoid $\mc G|_{C}$ by $s^{-1}(C)\cap r^{-1}(C)$, and note that if $C$ is invariant under $\mc G$ then this is a closed subgroupoid of $\mc G$, and that there is a natural restriction map from $C_{c}(\mc G)$ to $C_{c}(\mc G|_{C})$, which extends continuously to both maximal and regular representations.

Applying this construction to the groupoids attached to an inverse semigroup $S$, we obtain a natural surjective homomorphisms $C^{*}(G(S)) \rightarrow C^{*}(G_{tight}(S))$ and $C^{*}_r(G(S)) \rightarrow C^{*}_r(G_{tight}(S))$, because $\widehat{E}_{tight}$ is a closed invariant subset of $\widehat{E}$, for every inverse semigroup $S$ \cite{MR2419901}.

Finally, we can relate this construction directly to $C^{*}_{-}(S)$ defined in the previous section. By the work of Paterson \cite[Theorem 4.4.1, Theorem 4.4.2]{MR1724106} (see also Exel \cite{MR2419901} and Khoskham--Skandalis \cite{MR1900993}), we have the following isomorphisms:

\begin{equation*}
C^{*}(S) \cong C^{*}(G(S)),\quad C^{*}_r(S) \cong C^{*}_r(G(S)).
\end{equation*}

More is known concerning the connection between representations of $S$ and those of $G(S)$: in \cite{MR2419901}, Exel showed that representations of $S$ by partial isometries are in one-to-one correspondence with representations of $C_{c}(G(S))$, and that tight representations of $S$ are in one-to-one correspondence with representations of $C_{c}(G_{tight}(S))$, i.e those that factor through the canonical quotient map $C_{c}(G(S))\rightarrow C_{c}(G_{tight}(S))$. We will make use of this fact in the next section.

\subsection{Groupoid and Koopman representations of $F(\mc G)$}\label{Sect:groupoidrep}

We begin this section with a construction fundamental to later aspects of this paper.

Let $\mc G^{a}$ be the Boolean inverse semigroup attached to a ample \'etale groupoid $\mc G$. Then by appealing to Exel's reconstruction theorem \cite{MR2644910}, or one of the various non-commutative Stone dualities \cite{MR3077869}, the groupoid $G_{tight}(\mc G^{a})$ is topologically isomorphic to $\mc G$, and so the restriction of functions coupled with the identification of groupoid and semigroup \Cst-algebras described in the previous section induces a quotient homomorphism $\pi: \mathbb{C}\mc G^{a} \rightarrow C_{c}(\mc G)$, which maps a compact open bisection $U$ to its characteristic function $1_{U}\in C_{c}(\mc G)$.

\begin{Def}\label{Def:groupoidrep}
The canonical representation $\pi: F(\mc G) \rightarrow U(C_{c}(\mc G))$ obtained by restricting the construction above is called the \textit{canonical groupoid representation of $F(\mc G)$}.

A representation $\rho\colon F(\mc G)\to U(\mc H)$ is a \textit{groupoid representation of $F(\mc G)$} if it factors through the canonical groupoid representation $\pi$. If $\theta \colon C_c(\mc G)\to \mb B(\mc H)$ is any representation of $\mc G$, we denote the corresponding groupoid representation of $F(\mc G)$ by $\pi_\theta\coloneqq \theta\circ\pi$.
\end{Def}

% \begin{Rem}
% In the case that $\mc G$ is a minimal groupoid of germs and $A$ is the reduced \Cst-algebra of $\mc G$ this agrees with the representation defined following Takesaki \cite{MR1943007}, because, by Lemma \ref{Lem:germing2}, the orbits of $F(\mc G)$ acting on $\mc G^{(0)}$ agree with those of $\mc G$. 
% \end{Rem}

\begin{Rem}
This definition of a groupoid representation is nothing but a representation of $F(\mc G)$ that comes from a tight completion of Boolean inverse monoid $\mc G^{a}$. This follows from \cite[Theorem 13.3]{MR2419901} and the identifications from the previous section.
\end{Rem}

We are thus lead to consider quasi-regular representations:
\begin{equation*}
\rho_{x}: F(\mc G) \rightarrow \mb{B}(\ell^{2}(\mc G_{x}))
\end{equation*}
as they contain important information about the structure of the groupoid representation associated to the regular representation of $\mc G$.

These sorts of representation were studied independently by Birget \cite{MR2104771} and Nekrasheyvch \cite{MR2119267} in a particular example 
%(which will will study in detail in the next section) and appear in the work of \cite{}, and \cite{},
and was further studied by Artem Dudko and Rostislav Grigorchuk \cite{MR3627119} in the context of weakly branch groups.%\todo{add references}

Another key example of a representation is the \textit{Koopman representation} of $F(\mc G)$. Recall that a Borel measure $\mu$ on $\mc G^{(0)}$ is called quasi-invariant under $F(\mc G)$ if for every $g\in F(\mc G)$ the measures $g_*\mu$ and $\mu$ are equivalent. If $\mu$ be a $F(\mc G)$-quasi-invariant measure on $\mc G^{(0)}$, then we let $\kappa\colon F(\mc G)\to U(L^2(\mc G^{(0)},\mu))$ be the representation defined in the following way:
\begin{equation*}
\kappa(g)(f)(x) = \sqrt{\frac{d\mu(g^{-1}(x))}{d \mu (x)}}f(g^{-1}x),\quad f\in L^2(\mc G^{(0)},\mu)
\end{equation*}
This representation is extended to compact open bisections $U \in \mc G^{a}$, acting now by partial isometries rather than unitaries by the similar formula:
\begin{equation*}
\kappa(U)(f)(x) = \sqrt{\frac{d\mu(U^{*}(x))}{d \mu (x)}}f(U^{*}x).
\end{equation*}

\begin{Lemma}\label{Lem:Koopman}
The Koopman representation $\kappa : \mc G^{a} \rightarrow \mathbb{B}(L^{2}(\mc G^{(0)},\mu)$ defined above is a tight representation of $\mc G^{a}$, thus induces a groupoid representation of $F(\mc G)$ in the sense of Definition \ref{Def:groupoidrep}.
\end{Lemma}
\begin{proof}
As $E(\mc G^{a})$ is a Boolean algebra by \cite[Proposition 11.9]{MR2419901} it suffices to check that the representation induces a Boolean algebra homomorphism. This is true, since the idempotents are mapped precisely to characteristic functions of clopens via this representation.
\end{proof}

\subsection{Representations associated to subgroups of inverse semigroups}
For an inverse monoid $S$ there is no relationship between the reduced $C^{\ast}$-algebra of a inverse submonoid $T$ of $S$ and $C^{\ast}_{r}(S)$ in general. What we show in the following Proposition is that if $T$ is a group, then some connection exists. 

\begin{Prop}\label{prop:exact-subgroup-via-nuclearity}
Let $S$ be an inverse monoid and $U\subset S$ be a subgroup of $S$ with the same identity element. Then $C^*_r(U)$ is a subquotient of $C^*_r(S)$. 
\end{Prop}
\begin{proof}
 Let $\lambda\colon \mb C S\to \mb B(\ell^2 S)$ denote the left regular representation of $S$ and let
 \[
  A = \overline{\lambda(\mb C U)}\subset C^*_r(S)\subset \mb B(\ell^2 S).
 \] 
Consider the orthogonal projection $p\colon \ell^2 S\to \ell^2 U$ and observe that for every $g\in U$
\[
 p\lambda(g)p = \lambda_g\colon \ell^2U\to \ell^2 U,
\]
where $\lambda_g$ denotes the left regular representation of $U$ applied to $g\in U$. Moreover, for each $g \in U$, we know that $\ell^2 U\subset \ell^2 S$ is an invariant subspace for the unitary $\lambda(g)$, hence $p$ commutes with $\lambda(g)$ for all $g\in U$ and the map
\[
 q\colon A \to C^*_r(U),
\]
\[
 a\mapsto pap,
\]
is a $\ast$-homomorphism of $C^*$-algebras; it is surjective because it has dense range. Thus, $C^*_r(U)$ is a quotient of $A$. 
\end{proof}

\subsection{Induction of representations}
Denote by $\Res^{U}_{S}\pi$ the representation obtained by restricting $\pi$ to a group representation of the unital subgroup $U$.

\begin{Prop}\label{Prop:induction}
Let $S$ be an inverse monoid, and $U$ the group of units of $S$. For every unitary representation $\pi: U \rightarrow \mb{B}(\mc H)$ there is a representation $\Ind_{U}^{S}:S \rightarrow \mb{B}(\mc H')$ such that $\pi$ is contained in $\Res^{U}_{S}\Ind_{U}^{S}\pi$.
\end{Prop}
\begin{proof}
We construct a representation as follows. Observe that $U$ acts on $S$ by bijections, and thus the orbit space $S/U$ is well defined, let $p: S \rightarrow S/U$ be the map that assigns each $s \in S$ its $U$-orbit $sU$. We note that unlike the group case, these might not give bijective copies of $U$ -- denote by $U_{s}$ the stabiliser subgroup for each $s\in S$. Now consider the space of functions:
\begin{equation*}
\mathcal{F} := \lbrace f: S \rightarrow\mc H \mid f(su)=\pi(u)^{-1}f(s) \mbox{ and } p(\supp f) \mbox{ is finite} \rbrace. 
\end{equation*} 
As with group induction, our goal is to construct an inner product on this space. Let $f,f' \in \mathcal{F}$, and define:
\begin{equation*}
\langle f,f' \rangle := \sum_{x \in S/U} \langle f(x),f'(x) \rangle_{H}.
\end{equation*}
This is well defined, since the map
\begin{equation*}
x \mapsto \langle f(x), f'(x) \rangle
\end{equation*}
is constant on the right $U$-orbits. Take $\mc H'$ to the Hilbert space completion of $\mathcal{F}$ using this inner product, and define the map $\Ind_{U}^{S}\pi(s)$ using the left action of $S$ on itself:
\begin{equation*}
\Ind_{U}^{S}\pi(s)f(x) := f(s^{*}x).
\end{equation*}
This is clearly a homomorphism of $S$, and each $\Ind_{U}^{S}\pi(s)$ is a partial isometry as it defines a linear map between the subspaces: $\mc H'_{s^{*}s}$ and $\mc H'_{ss^{*}}$, where $\mc H'_{e}$ is the closure in $\mc H'$ of all the functions containing the orbit of $e$ in their support.

The final claim concerning the restriction $\Res^{U}_{S}\Ind_{U}^{S}\pi : U \rightarrow U(\mc H')$ follows by restricting the representation to the $U$-invariant subspace $\mc H'_{1_{S}}$, this has a free $U$-action, in particular this gives rise to a multiple of original representation $\pi$.

\end{proof}

In the case of the representation $\lambda_{S}$ we can describe the weak equivalence class of $\pi=\Res^{U}_{S}\lambda_{S}$.

\begin{Prop}\label{Prop:leftregular}
$\pi=\Res^{U}_{S}\lambda_{S}$ is weakly equivalent to $\bigoplus_{e \in E(S)/U} \lambda_{U/U_{e}}$, where $\lambda_{U/U_{e}}$ is the quasiregular representation of $U$ on $\ell^{2}(U/U_{e})$.
\end{Prop}
\begin{proof}
Splitting $S$ into right $U$-orbits we decompose $\ell^{2}(S)$ as $\bigoplus_{s \in S/U}\ell^{2}(sU)$. Letting $U_{s}$ be the stabiliser of $sU$ under the left action, we obtain that the restriction of $\lambda_{S}$ decomposes as $\bigoplus_{s \in S/U} \lambda_{U/U_{s}}$, however, we note that since $U_{s}=U_{ss^{*}}$, that this summand is weakly equivalent to the one in the claim.
\end{proof}

\section{Comparing some natural representations of the full group}
The goal of this section is to relate these representations of $F(\mc G)$ following the ideas of \cite{Grig-Dudko} and \cite{MR3627119} where the corresponding theory was developed in the weakly branch case. The main result is the following:

\begin{Thm}\label{Thm:irred}
Let $\mc G$ be an ample effective groupoid, $F(\mc G)$ its topological full group and $\lambda$ the left regular representation of $\mc G$. Then for the representations $\pi_{\lambda}, \pi_{\kappa}$ and $\rho_{x}$ defined above we have:
\begin{enumerate}
\item If $\mc G$ is minimal, then: \begin{enumerate} 
\item for every $x \in X$ we have $\pi_{\lambda} \sim \rho_{x}$;
\item if $\mc G$ has infinite orbits, then the representation $\rho_{x}$ is irreducible for each $x \in X$;
\item for $x,y \in X$, with $x$ and $y$ in distinct $\mc G$-orbits, $\rho_{x}$ and $\rho_{y}$ are not unitarily equivalent. 
\end{enumerate}
\item If $\mc G$ is topologically amenable, then $\pi_{\kappa} \prec \pi_{\lambda}$.
\item If $\mc G$ is minimal and topologically amenable, then $\pi_{\kappa} \sim \pi_{\lambda}$.
\end{enumerate}
\end{Thm}
\begin{proof}[Proof of (i) (a), (ii) and (iii)]
For (i) (a), we observe that $\Vert \pi(g) \Vert_{r}$ is the reduced groupoid norm of the element $g \in C_{c}(\mc G)$, and thus it bounds above the norm $\Vert \lambda_{x}(\pi(g)) \Vert$ for each $x$. Note that the representation $\lambda_{x}\circ \pi$ agrees with the quasi-regular representation $\rho_{x}$ of $F(\mc G)$. By a result of Khoshkam and Skandalis \cite[Corollary 2.4]{MR1900993}, it is sufficient to compute the supremum  over a dense subset of the unit space, such as any single $\mc G$-orbit. But by Proposition \ref{Lem:germing2}, these are precisely the $F(\mc G)$-orbits.

For (ii), we observe that as $\mc G$ is amenable, the reduced \Cst-algebra $C^{*}_{r}(\mc G)$ is isomorphic to the maximal \Cst-algebra $C^{*}(\mc G)$, and so it satisfies the universal property of the maximal completion. By Lemma \ref{Lem:Koopman}, the representation $\kappa$ is a groupoid representation of $F(\mc G)$, thus it is obtained through some completion of $C_{c}(\mc G)$ \cite[Theorem 13.3]{MR2419901}. Then $\Vert \kappa(g) \Vert = \Vert \pi(g) \Vert_{\mb{B}(L^{2}(X,\mu))} \leqslant \Vert \pi(g) \Vert_{\lambda} = \Vert \pi_{\lambda}(g) \Vert$, where the middle inequality follows from the fact that the reduced norm and maximal norm on $C_{c}(\mc G)$ agree. Thus $\pi_{\kappa} \prec \pi_{\lambda}$.

(iii) follows from (ii) and the fact that $C^{*}_{r}(\mc G)$ is a simple \Cst-algebra when $\mc G$ is minimal.%\todo{add reference}. 
\end{proof}

To show (i) (b) and (c) we use Mackey's criterion for irreducibility and disjointness of quasi-regular representations in terms of commensurators of subgroups and the notion of quasi-conjugacy.

\begin{Def}\label{Def:commensurator}
Recall that $H,K < G$ are \textit{commensurable} if $H\cap K$ has finite index in $H$ and $K$ and are \textit{quasi-conjugate} if $H^{g}:=gHg^{-1}$ is commensurate with $K$ for some $g \in G$. the commensurator is defined, for $H \leqslant G$, as:
\begin{equation*}
\Comm_{H}(G) := \lbrace g \in G \mid H\cap H^{g} \mbox{ has finite index in } H \mbox{ and } H^{g}\rbrace.
\end{equation*}
\end{Def}

Let us recall Mackey's criterion:
\begin{Thm}[{\cite[Corollary 7]{MR0396826}}]
Let $G$ be a countable discrete group, $H$, $K$ be subgroups of $G$.
\begin{enumerate}
\item The quasi-regular representation $\rho_{G/H}$ is irreducible if and only if $\Comm_{G}(H)=H$.
\item The quasi-regular representations $\rho_{G/H}$ and $\rho_{G/K}$ are unitarily equivalent if and only if $H$ and $K$ are quasi-conjugate in $G$.
\end{enumerate}\qed
\end{Thm} 

\begin{Def}
Let $U\subseteq \mc G^{(0)}$ be clopen. The rigid stabiliser $F(\mc G)_{(U)}$ of $U$ is the subgroup of $F(\mc G)$ consisting of all elements which act trivially on $\mc G^{(0)}\setminus U$.
\end{Def}

The following technical lemma will be used in the sequel.

\begin{Lemma}\label{Lem:technicalminimal}
Let $\mc G$ be a minimal effective groupoid with infinite orbits. Then
\begin{enumerate}
\item each rigid stabiliser $F(\mc G)_{(U)}$ is isomorphic to the full group $F(\mc G|_U)$ of the restriction of $\mc G$ to $U$.
% \item for every pair of clopen subsets $U, V \not = \mc G^{(0)}$, there is an element $g \in F(\mc G)$ such that $g(U)\cap V \not = \varnothing$.
\item for every pair of non-empty clopen subsets $U,V \not = \mc G^{(0)}$ and every point $x \in V$ there exists a $g \in F(\mc G)$ such that $x\in g(U)\cap V$ (in particular, $g(U)\cap V\neq\varnothing$).
\end{enumerate}
\end{Lemma}
\begin{proof}
To show (i), we observe that for every non-empty clopen $U\subset \mc G^{(0)}$, the restriction groupoid $\mc G|_{U}$ is also effective, and the rigid stabilisers $F(\mc G)_{(U)}$ are precisely the groups $F(\mc G|_{U})$ extended trivially on the complement of $U$. By \cite[Lemma 3.1]{MR1841750}, these contain arbitrarily long finite orbits, thus $F(\mc G)_{(U)}$ is infinite.

For (ii): Notice that due to the minimality of the groupoid, for every non-empty clopen $U$, the union $\bigcup_{g \in F(\mc G)}g(U)$ is equal to $\mc G^{(0)}$, thus one of them contains $x$. This proves (ii).
\end{proof}

The following Lemmas are modifications of the work of Bartholdi and Grigorchuk \cite{MR1841750} as well as Dudko and Grigochuk \cite{Grig-Dudko} into the setting of full groups.

\begin{Lemma}\label{Lem:irred}
Let $\mc G$ be a minimal effective groupoid with infinite orbits, $\Gamma=F(\mc G)$ its topological full group, $x \in \mc G^{(0)}$ and $S_{x}=\Stab(x)<\Gamma$. Then $\Comm_{\Gamma}(S_{x})=S_{x}$.
\end{Lemma}
\begin{proof}
Let $g \in \Gamma\setminus S_{x}$. We will show the subgroup $S_{x}^{g}\cap S_{x}$ has infinite index in $S^{g}_{x}$ by showing that the orbit $S_x^g\cdot x$ is infinite.

First we observe that as the base space is Hausdorff and totally disconnected, there is a clopen neighbourhood $U$ of $x$ such that $g(U)\cap U = \varnothing$. Now $(S_{x}^{g})_{(U)} = F(\mc G)_{(U)}$ since $g(U)\subset \mc G^{(0)}\setminus U$ by assumption, wherefore any element of $F(\mc G)$ that fixes $\mc G^{(0)}\setminus U$ pointwise must fix $g(U)$ pointwise, thus belonging to $S_{gx} = S_x^g$. Hence we have $(S_{x}\cap S_{x}^{g})_{(U)} = S_{x}\cap F(\mc G)_{(U)}= (S_{x})_{(U)}$.

The orbit of $F(\mc G)_{(U)} \cdot x$ is infinite by Lemma \ref{Lem:technicalminimal} (i) and minimality of the restriction groupoid. As moreover $F(\mc G)_{(U)} \leqslant S_x^g$, the orbit $S_x^g\cdot x$ also follows infinite, which concludes the proof.
\end{proof}

\begin{Lemma}\label{Lem:quasiconjugate}
If $x,y \in X$ belong to distinct $\mc G$-orbits, then $S_{x}$ and $S_{y}$ are not quasi-conjugate in $G$.
\end{Lemma}
\begin{proof}
Our goal is to show that $S_{x}y$ and $S_{y}x$ are infinite, which is equivalent to quasi-conjugacy after replacing $x$ by $gx$ for some appropriate $g$. This follows from the definition \ref{Def:commensurator}.

Let $\lbrace W_{i} \rbrace_{i\in \mathbb{N}}$ be a nested family of neighbourhoods of $x$ with $\cap_{i} W_{i} = \lbrace x \rbrace$. Then, pick an $i$ large enough such that $y \not \in W_{i}$, and find a clopen neighbourhood $V\ni y$ such that $V \cap W_{i}=\emptyset$ (which can be done because $X$ is regular). 

Now we will make repeated use of Lemma \ref{Lem:technicalminimal} (ii) to complete the proof. Let $h \in F(\mc G)_{W_{i}}$ be some element that is not trivial, and then for each $j\geq i$, find a $k_{j}$ such that $k_{j}(\supp(h))\cap W_{j}$ contains $x$, and is thus not empty. Let $V_{j}$ denote $k_{j}^{-1}(k_{j}(\supp(h)\cap W_{j})$.

Set, for each $j$, $Z_{j}=W_{j+1}\setminus W_{j}$, which is clopen as all the $W_{j}$ are. Then again, apply Lemma \ref{Lem:technicalminimal} (ii) to obtain a $z_{j}$ such that $z_{j}(Z_{j})\cap h(V_{j})$ is not empty and contains $hk_{j}^{-1}x$. Then the elements $z_{j}^{-1}hk_{j}^{-1}$ does not stabilise $x$ (so it is not trivial), belongs to $F(\mc G)_{W_{i}}$ (and thus $S_{y}$, since $y \not \in W_{i}$) and maps $x$ into $Z_{j}$.

Note that for each $j,l \in \mathbb{N}$, the sets $Z_{j}$ and $Z_{l}$ have empty intersection, so the points $z_{j}^{-1}hk_{j}^{-1}x$ are all distinct, and this provides infinitely many distinct points in the orbit $S_{y}x$. By symmetry, we can construct an infinite subset of the orbit $S_{x}y$.
\end{proof}

Now, parts (b) and (c) of Theorem \ref{Thm:irred} (i) follow from Mackey's criterion directly.

% \newpage

\section{Amenability and \Cst-simplicity of $F(\mc G)$}

\subsection{Invariant means and amenability}

\begin{Def}
Let $S$ be a Boolean inverse semigroup (such as $\mc G^{a}$). Then $S$ has an \textit{invariant mean} if there is a function $\mu: E(S) \rightarrow [0,\infty)$ such that:
\begin{enumerate}
\item for any $s \in S$, we have that $\mu(s^{*}s)=\mu(ss^{*})$,
\item if $e$ and $f$ are orthogonal idempotents, then $\mu(e \vee f) = \mu(e)+\mu(f)$.
\end{enumerate}
A mean $\mu$ is \textit{normalised} at $e \in E(S)$ if $\mu(e)=1$. By convention, if we do not mention any normalisation then we suppose $\mu(1)=1$, and \textit{faithful} if $\mu(e)=0$ implies $e=0$.
\end{Def}

The existence of invariant means on Boolean inverse monoids is investigated by Kudryavtseva, Lawson, Lenz, and Resende \cite{MR3448402}, and we will go into detail into the ideas in that work in the next section. Additionally in the context of an ample groupoid $\mc G$, Starling \cite{MR3548134} addressed the existence of invariant means on $\mc G^{a}$ in the context of both invariant measures on the base space of $\mc G$ and traces on the reduced $C^{*}$-algebra $C^{*}_{r}(\mc G)$.

% \begin{Thm}
% Suppose $\mc G$ is an ample Hausdorff principal groupoid, then the following sets are in bijection:
% \begin{itemize}
% \item $\mathrm{IM}(\mc G):=\lbrace \mu \in \Prob(G^{(0)}) \mid \mu(s(U))=\mu(r(U))\text{ for all } U \in \mc G^{a}\rbrace$
% \item $\mathrm{M}(\mc G^{a}):=\lbrace m \mid m \text{ a mean on } E(\mc G^{a})\rbrace$.
% \item $\mathrm{T}(C^{*}_{r}(\mc G)) = \lbrace \tau\colon C^{*}_{r}(\mc G)\to \mb C\mid \tau\text{ a tracial state}\rbrace$.
% \end{itemize}
% \end{Thm}

% However, asking that the groupoid $\mc G$ preserves a measure on its own base space $\mc G^{(0)}$ is not the only possible notion of amenability one could consider - the other is \textit{topological amenability}, which we recall below:

% \begin{Def}(topological amenability)

% \end{Def}

We now relate the existence of an invariant mean on $\mc G^{a}$ to the amenability of $F(\mc G)$ in the minimal case, partly addressing Remark 2.25 in \cite{MR3448402}.

\begin{Thm}\label{Thm:amenable}
Let $\mc G$ be a minimal ample effective groupoid. Then $F(\mc G)$ is amenable if and only if the following three conditions are satisfied:
\begin{enumerate}
\item $\mc G^{a}$ admits an invariant mean,
\item $\mc G$ is topologically amenable and
\item the rigid stabiliser $F(\mc G)_{(U)}$ is amenable for each clopen subset $U \subset \mc G^{(0)}$.
\end{enumerate}

\end{Thm}
\begin{proof}
Suppose $F(\mc G)$ is amenable, then it is immediate that (iii) holds, since subgroups of amenable groups are amenable. To conclude (ii), we appeal to minimality: as $\mc G$ minimal, we can apply Proposition \ref{Lem:germing2} to conclude that $\mc G \cong \Germs(F(\mc G)\curvearrowright \mc G^{(0)})$ is a groupoid quotient of $F(\mc G)\curvearrowright \mc G^{(0)}$. Since $F(\mc G)$ is amenable, the transformation groupoid $F(\mc G)\curvearrowright \mc G^{(0)}$ is topologically amenable. By applying Proposition 5.1.2 from \cite{MR1799683}, which states that topological amenability is groupoid extension closed, we can conclude that $\mc G$ is topologically amenable. Finally, (i) is Proposition 2.24 in \cite{MR3448402} (which was the original motivation for considering this problem).

We now show the other direction by supposing conditions (i), (ii) and (iii). We will appeal to representation theory to conclude the result. We analyse groupoid representations of $F(\mc G)$: the first step is to use group to inverse semigroup induction to understand the natural groupoid representation in $C^{*}_{r}(\mc G)$, and the second step involves studying the Koopman representation of $F(\mc G)$. We combine these using Theorem \ref{Thm:irred} using assumption (ii).

Let $\pi:=\Res_{F(\mc G)}^{\mc G^{a}}\lambda_{\mc G^{a}}$ be the representation of $F(\mc G)$ obtained by completing $\mb C F(\mc G)$ in $C^{*}_{r}(\mc G^{a})$. By Proposition \ref{Prop:leftregular}, the representation $\tilde{\pi}$ is weakly equivalent to $\bigoplus_{U} \lambda_{F(\mc G)/F(\mc G)_{(U)}}$, where $F(\mc G)_{(U)}$ is the rigid stabiliser of $U$. However, since each $F(\mc G)_{(U)}$ is amenable, each representation $\lambda_{F(\mc G)/F(\mc G)_{(U)}}$ is weakly contained in the left regular representation $\lambda_{F(\mc G)}$ (by induction), so we know that $\pi$ is weakly equivalent to $\lambda_{F(\mc G)}$.

Consider the groupoid representation $\pi_{\lambda}$ of $F(\mc G)$ in $C^{*}_{r}(\mc G)$. This representation $\pi_{\lambda}$ is weakly contained in $\pi$, since the algebraic map:
\begin{equation*}
\mb C \mc G^{a} \rightarrow C_{c}(\mc G)
\end{equation*}
extends surjectively to reduced completions, so for any $g\in F(\mc G)$, there is an inequality of norms $\Vert \pi_{\lambda} (g) \Vert \leqslant \Vert \pi(g) \Vert$. Appealing to a result of Dixmier \cite[Theorem F.4.4]{MR2415834}, this is equivalent to weak containment of $\pi_{\lambda}$ in $\pi$.

The second step is to conclude something about the Koopman representation of $F(\mc G)$, which is also a groupoid representation by Lemma \ref{Lem:Koopman}. As $\mc G^{a}$ admits an invariant mean $\mc G$ admits an invariant measure by Corollary 4.12 in \cite{MR3548134}. Note that this measure is automatically $F(\mc G)$-invariant, and so by Theorem 5.7 in \cite{Bekka}, the representation $\pi_{\kappa}$ has almost invariant vectors -- that is, $1_{F(\mc G)}$ is weakly contained in $\pi_{\kappa}$. 

Since $\mc G$ is topologically amenable, we can apply part (iii) of Theorem \ref{Thm:irred} to obtain the following chain of weak containments and equivalences:
\begin{equation*}
1_{F(\mc G)} \prec \pi_{\kappa} \sim \pi_{\lambda} \prec \pi \sim \lambda_{F(\mc G)}.
\end{equation*}
Thus $F(\mc G)$ is amenable (since weak containment of the trivial representation in the left regular representation is a characterisation).
\end{proof}

We remark that it is possible for $F(\mc G)$ to be non-amenable, whilst $\mc G^{a}$ admits an invariant mean -- due to results of Starling \cite{MR3548134}, Kerr and Nowak \cite{MR2974211}, it is the case for any residually finite action of a free group (or more generally, an action of a non-amenable group preserving a probability measure on the Cantor set).

This also naturally connects with the work of Uffe Haagerup and Kristian Olesen \cite{Haagerup-Olesen} part of which is concerned with \Cst-simplicty of Thompson's group $T$ and the amenability of $F$, as well as the work of Le Boudec and Matte-Bon \cite{LB-MB}.

\begin{Lemma}\label{Lem:algkern}
Let $\mc G$ be a minimal ample \'etale groupoid. Then the representation $\pi:G \rightarrow C_{c}(\mc G)$ does not extend faithfully to $\mathbb{C}G$ for any subgroup $G$ of $F(\mc G)$ with non-trivial rigid stabilisers. Thus it is not injective for any faithful groupoid representation of $G$.
\end{Lemma}
\begin{proof}
Note that by minimality, for every clopen subset $U \subset X$, the subgroup $F(\mc G|_{U})$ is infinite. Thus, we can partition of $X$ into $U$ and $U^{c}$, and then consider $g_{1} \in F(\mc G|_{U})$ and $g_{2} \in F(\mc G|_{U^{c}})$. These elements are faithfully represented in $C_{c}(\mc G)$ via the construction at the beginning of Section \ref{Sect:groupoidrep}, and satisfy the formula:
\begin{equation*}
1+g_{1}g_{2}=g_{1}+g_{2}.
\end{equation*}
Thus, the element $a = g_{1}g_{2}-g_{1}-g_{2}+1=0$ in any faithful groupoid representation of $F(\mc G)$ extended linearly to $\mb C F(\mc G)$.
\end{proof}

% The subgroup $F(\mc G|_{U})$ is the \textit{rigid stabiliser} of $U$ in the context of \cite{}.

The following theorem can be regarded as a ``groupoid version'' of the ideas of Haagerup--Oleson \cite{Haagerup-Olesen}, Birget \cite{MR2104771} and Nekrashevych \cite{MR2119267}. It also connects with recent work of Le Boudec and Matte Bon \cite{LB-MB}.

\begin{Thm}\label{Thm:not-c-star-simple}
Let $\mc G$ be a minimal effective groupoid. If the rigid stabiliser $F(\mc G)_{(U)}$ is amenable for every clopen subset $U$, then the group $F(\mc G)$ is not $C^{*}$-simple.
\end{Thm}
\begin{proof}
This uses part of the argument for Theorem \ref{Thm:amenable}, but we repeat it here. Let $\pi:=\Res_{F(\mc G)}^{\mc G^{a}}\lambda_{\mc G^{a}}$ be the representation of $F(\mc G)$ obtained by completing $\mb C F(\mc G)$ in $C^{*}_{r}(\mc G^{a})$. By Proposition \ref{Prop:leftregular} $\tilde{\pi}$ is weakly equivalent to $\bigoplus_{U} \lambda_{F(\mc G)/F(\mc G)_{(U)}}$. However, since each $F(\mc G)_{(U)}$ is amenable, each representation $\lambda_{F(\mc G)/F(\mc G)_{(U)}}$ is weakly contained in the left regular representation $\lambda_{F(\mc G)}$ (by induction), so we know that the groupoid representation $\pi_{\lambda}$ is weakly contained in $\lambda_{F(\mc G)}$. However, the representation $\pi_{\lambda}$ has an algebraic kernel by Lemma \ref{Lem:algkern}, so $C^{*}_{\pi_{\lambda}}(F(\mc G)$ is a proper quotient of $C^{*}_{r}(F(\mc G))$, which completes the proof. 
\end{proof}

% Note that this is a conceptualisation of the ideas of Haagerup--Oleson \cite{Haagerup-Oleson}, which appears also in Nekrasheyvch \cite{MR2119267}, Birget \cite{MR2104771} and in other work building on that such as \cite{} and \cite{} - one can also find references pointing to Lemma \ref{Lem:algkern} in the work of Grigorchuk--Musat--Rørdam \cite{GMR} on just infinite $C^{*}$-algebras.

\section{A comparison of type semigroups}

% \subsection{A brief history of the type semigroup}This points a lot to WAGON, and states the original theorem of Tarski - then points to the modifications in recent 
% years - such as looking at C-S--dlH--Grig.

In this section, we combine the results of \cite{MR3448402},  \cite{MR3548134} and \cite{Bonicke-Li} to obtain a unified Tarski alternative for ample groupoids both in terms of the associated semigroups and \Cst-algebras. %We refer to the above articles for the relevant definitions.

\subsection{A plethora of definitions}
First, we recall some preliminaries concerning ordered commutative monoids, and then give two different definitions of the type semigroup for an \'etale groupoid, and show they agree. For more information and further relevant notions see \cite{MR3448402}.

\begin{Def}
Let $M$ be a commutative monoid with addition $+$ and identity $0$. Let $a, b\in M$. Then:
\begin{itemize}
\item the algebraic preorder is defined as $a \leqslant b$ if and only if there exists a $c \in M$ such that $a+c =b$. 
\item we write $na$ to represent the sum $\overbrace{a+...+a}^n$.
\item $u \in M$ is an order unit if for every $a\in M$ there is an $n$ such that $a\leqslant nu$.
\item Note that $[0,\infty)$ is an commutative monoid with order unit $1$.
\end{itemize}
\end{Def}

We begin by giving the ad-hoc definition, following \cite{Bonicke-Li}.% and \cite{}.

Consider the set
\begin{equation*}
\mathcal{A} := \lbrace \cup_{i=1}^{n} A_{i}\times \lbrace i\rbrace \mid n \in \mathbb{N}, A_{i}\in G^{a}, A_{i}\subset \mc G^{(0)} \rbrace 
\end{equation*}
which collects up, using multiplicity, clopen subsets of the unit space of $\mc G$. We now construct an equivalence relation $\sim$ on $\mc A$ as follows. Given two sets $A=\bigcup_{i=1}^{n} A_{i}\times \lbrace i \rbrace$ and $B = \bigcup_{j=1}^{m}B_{j}\times \lbrace j \rbrace$, we write $A\sim B$ if there exists an $l \in \mathbb{N}$ and $\mc G$-bisections $W_{1},...,W_{l}\in \mc G^{a}$ and numbers $n_{1},...,n_{l},m_{1},...,m_{l} \in \mathbb{N}$ such that:
\begin{equation*}
A = \bigsqcup_{k=1,...,l} s(W_{k})\times {n_{k}} \mbox{ and } B = \bigsqcup_{k=1}^{l}r(W_{k})\times \lbrace m_{k} \rbrace.
\end{equation*}
This is an equivalence relation (Lemma 5.2 in \cite{Bonicke-Li}) and leads to the following definition:
\begin{Def}
The \textit{type semigroup} of $\mc G$, denoted $S(\mc G, \mc G^{a})$ is the set $\mathcal{A}/\sim$ with the product: 
\begin{equation*}
\left[\bigcup_{i=1}^{n} A_{i}\times \lbrace i \rbrace\right]+\left[\bigcup_{j=1}^{m} B_{j}\times \lbrace j \rbrace\right] = \left[\bigcup_{i=1}^{n} A_{i}\times \lbrace i \rbrace \cup \bigcup_{j=1}^{m} B_{j}\times \lbrace n+j \rbrace\right]
\end{equation*}
where $[A]$ is the equivalence class of $A$ in $S(\mc G, \mc G^{a})$.
\end{Def}

% \subsubsection{A second operator algebraic definition}

\subsubsection{The Boolean inverse monoid definition}

Secondly, we give the definition for Boolean inverse monoids originally considered by Kudryatseva--Lawson--Lenz--Resende in \cite{MR3448402}. To do this, we must must first recall two of \textit{Green's relations} in the context of inverse semigroups:

\begin{Def}[Relations $\mc D$ and $\mc J$]
Let $S$ be an inverse semigroup and let $e,f \in E(S)$ be its idempotents. Then:
\begin{itemize}
\item $e \mc D f$ if there is an $s \in S$ such that $e=s^{*}s$ and $f=ss^{*}$,
\item $e \leq_{\mc J} f$ if there exists $e'\in E(S), e' \leqslant f$ such that $e\mc D e'$ and
\item $e \mc J f$ if $e \leq_{J} f$ and $f \leq_{J} e$. 
\end{itemize}
\end{Def}
Note that in general, $\mc D \subset \mc J$ as relations on $E(S)$, and the equality $\mc D = \mc J$ is an inverse semigroup expression of the Sch\"{o}der--Bernstein theorem.

%(this is given as a remark on page 3 of \cite{}) - to see why, let $X$ be a countable set. Consider partial bijections $I(X)$ associated to subsets $A,B \subset X$, then $1_{A}\leq_{J} 1_{B}$ in the power set of $X$ if there is an injection from $A$ to $B$, and $1_{A} \mc D 1_{B}$ if there is a bijection between $A$ and $B$. 

% In \cite{}. the authors considered the quotient $E(S)/\mc D$ when $S$ is a Boolean inverse monoid - this may not be a commutative monoid however - the issue is the lack of stabilisation: there could be elements $e,f \in E(S)$ such that there are no pair of idempotents $e', f'$ that are below $e$ and $f$ respectively such that $e' \perp f'$. When such a pair exists, we say $E(S)$ is orthogonally separating, and can define a monoid operation $[e]\oplus [f] = [e'\vee f']$ for the corresponding orthogonal pair - note that this is precisely the operation in operator K-theory for idempotents.

% Let's note the following lemma; the proof is straightforward.
% \begin{Lemma}
% Let $e,f \in E(S)$ be idempotents. Then $e\sim f$ in any \Cst-algebra completion if $[e]=[f]$.
% \end{Lemma}

% Not all Boolean inverse monoids are idempotent separating (in the same way not all \Cst-algebras are stable), but we can define a stable completion in which $S$ embeds in the following way. 

\begin{Def}
A $n\times m$ \textit{generalised rook matrix} over $S$ is a $n\times m$ matrix with values in $S$ that satisfies:
\begin{enumerate}
\item if $a, b \in S$ are in distinct columns, then $a^{*}b=0$, i.e $aa^{*}\perp bb^{*}$.
\item if $a ,b \in S$ are in distinct rows, then $ab^{*}=0$, i.e $a^{*}a \perp b^{*}b$.
\item in the case that both $m$ and $n$ are infinite, we require that there are at most finitely many non-zero entries.
\end{enumerate} 
\end{Def}

When $n$ is finite, let $I_{n}$ be the identity matrix, i.e the matrix with all entries on the diagonal the unit in the Boolean inverse monoid $S$. Note that for an $m\times n$ matrix $A$ and a $n\times p $ matrix $B$, the product defined by:
\begin{equation*}
AB_{ij}=\bigvee_{k}a_{ik}b_{kj}
\end{equation*}
is well defined, and a rook matrix. Proposition 3.5 in \cite{MR3448402} shows that the collection: $M_{n}(S)$ (and $M_{\omega}(S)$) of finite, $n\times n$  (or infinite, with finitely many non-zero entry) rook matrices is itself a Boolean inverse monoid (resp. semigroup) when $S$ is.

Note that the map $\Delta(a_{1},...,a_{n}) = \diag(a_{1},...,a_{n})$ defines an embedding of $S^{n}$ into $M_{n}(S)$. Additionally, we extend this mapping to $\Delta_{\omega}(a_{1},...,a_{n})$ by defining the image to be the matrix with the first entries on its diagonal to be $a_{1},...,a_{n}$ and all the other entries $0$.

\begin{Lemma}[Lemmas 3.8 and 3.9 from \cite{MR3448402}]
If $S$ is a Boolean inverse monoid, then
\begin{itemize}
\item $S$ is isomorphic to the local monoid $\Delta_{\omega}(1)M_{\omega}(S)\Delta_{\omega}(1)$.
\item $S$ has an invariant mean if and only if $M_{\omega}(S)$ has an invariant mean that normalises $\Delta_{\omega}(1)$.
\item $M_{\omega}(S)$ is idempotent separating.
\end{itemize}\qed
\end{Lemma} 

\begin{Def}
The type monoid of $S$, denoted $T(S)$, is the monoid $E(M_{\omega}(S))/\mc D$, equipped with the operation defined above. Define a map $\delta: E(S) \rightarrow T(S)$ by $\delta(e)=[\Delta_{\omega}(e)]$, and let $\textbf{u}:=\delta(1)$.
\end{Def}

\subsection{The comparison theorem}

Our first result in this section is the following comparison theorem.

\begin{Thm}
Let $\mc G$ be an ample, \'etale groupoid. Then $S(\mc G,\mc G^{a}) \cong T(\mc G^{a}) $ as ordered commutative monoids.
\end{Thm}
\begin{proof}
Define a map $\mc A \rightarrow M_{\omega}(\mc G^{a})$ as follows:
\begin{equation*}
\bigcup_{i=1}^{n} A_{i}\times \lbrace i \rbrace \mapsto \Delta_{\omega}(A_{1},...,A_{n}),
\end{equation*}
which is well defined since $A_{i} \in \mc G^{a}$ and are idempotents there.% This map is a bijection...

It suffices to check that the relation defined on $\mc A$ is precisely the $\mc D$-relation on $E(M_{\omega}(\mc G^{a}))$. Note that if $A=\bigcup_{i=1}^{n} A_{i}\times \lbrace i \rbrace$ and $B = \bigcup_{j=1}^{m}B_{j}\times \lbrace j \rbrace$ and $A\sim B$, then there exists an $l \in \mathbb{N}$ and $\mc G$-bisections $W_{1},...,W_{l}\in \mc G^{a}$ and numbers $n_{1},...,n_{l},m_{1},...,m_{l} \in \mathbb{N}$ such that
\begin{equation*}
A = \bigsqcup_{k=1,...,l} s(W_{k})\times {n_{k}} \mbox{ and } B = \bigsqcup_{k=1}^{l}r(W_{k})\times \lbrace m_{k} \rbrace.
\end{equation*}
Under the map above, $A$ maps to $\Delta_{\omega}(e_{1},....,.e_{h})$, and $B$ maps to $\Delta_{\omega}(f_{1},....,.f_{h})$, where $h \in \mb N$ is the largest of the numbers $n_{1},...,n_{l},m_{1},...,m_{l}$ and $e_{i}$ (resp. $f_{i}$) is the subset $s(W_{i})$ (resp. $r(W_{i})$) if $i \in \lbrace n_{1},...,n_{l} \rbrace$ (resp. $i \in \lbrace m_{1},...,m_{l} \rbrace$) and $0 \in E(\mc G^{a})$ otherwise. Note that after reshuffling on the diagonal (something obtained using row or column operations, which are $\mc D$-equivalent) the matrix $\diag(W_{1},...,W_{k}) \in M_{\omega}(\mc G^{a})$ then implements the $\mc D$-relation. This shows that the map is well defined at the level of type semigroups.
\end{proof}

% This result lets us utilise techniques from the general theory of Boolean inverse monoids and their $C^{*}$-algebras developed over the papers \cite{}, ... , \cite{} and obtain results in the context useful to operator algebraists - this is the goal of the next section.

Next, we establish various forms of paradoxicality and Tarski alternatives.

% \begin{Def}
% Tarski matrices over a Boolean inverse monoid, paradoxical sets for groupoids
% \end{Def}

% \begin{Lemma}\label{Lem:paradoxical}
% Combining Bonicke--Li and Kudryatseva--Lawson--Lenz--Resende to get a Tarski matrix if and only if type monoid result if and only if paradoxical set.
% \end{Lemma}

% The upside of giving this ``stabilised" definition of \cite{} is the following observation, and Proposition (which is already in the \Cst-algebra literature, c.f Lemma x.xx in \cite{})

\begin{Prop}
Let $\mc G$ be a Hausdorff ample \'etale groupoid and let $\pi$ be the algebraic groupoid representation of $\mc G^{a}$ in $C_{c}(\mc G^{a})$. Suppose there is a Tarski matrix $A$ of degree $k$ over $\mc G^{a}$, then $AA^*$ is an infinite projection in $M_{k}(A)$ for any completion $A$ of $C_{c}(\mc G)$.  
\end{Prop}
\begin{proof}
Amplify the map $\pi$ to $\pi^{(k)}:M_{k}(\mc G^{a}) \rightarrow M_{k}(C_{c}(\mc G))$-entrywise, and equip $M_{k}(A)$ with the matrix norm that makes $M_{k}(C_{c}(\mc G))$ dense there. Then $\pi^{(k+1)}(A^*A) = 1_{k+1}$ and so $\pi^{(k+1)}(A^*A)$ is a proper subprojection (it lies in $M_k(A)$) equivalent to $1_{k+1}$.
\end{proof}

\begin{Cor}
Let $\mc G$ be an ample, minimal, Hausdorff \'etale groupoid with compact unit space. If there is a Tarski matrix of any degree over $\mc G^{a}$, then no $C^{*}$-completion of $C_{c}(\mc G)$ is stably finite.
\end{Cor}
\begin{proof}
By minimality, if there is a Tarski matrix of any degree, it can be taken to have full support, i.e there is a Tarski matrix $A$ of some degree with $A^{*}A=1$. From the proposition we get that $1$ would be infinite in $M_{k}(A)$ for any $A$ is any $C^{*}$-completion of $C_{c}(\mc G)$, since it is algebraically infinite. This contradicts stable finiteness of $C^{*}_{r}(\mc G)$.
\end{proof}

% The following is a modification and improvement of the classical theorem of Tarski \cite{}, which is due to Kudryatseva--Lawson--Lenz--Resende.

% \begin{Thm}

% \end{Thm}

We can now collect and unify various ``Tarski alternative'' type results both from $\Cst$-algebras and inverse semigroup theory \cite{MR3448402,MR3548134,Bonicke-Li,Rainone-Sims}.

\begin{Thm}[Tarski alternatives of \cite{MR3448402},  \cite{MR3548134} and \cite{Bonicke-Li}, combined]
Let $\mc G$ be an ample, minimal, Hausdorff groupoid with compact unit space. Then the following are equivalent:
\begin{enumerate}
\item $C^{*}_{r}(\mc G)$ admits a faithful tracial state,
\item $C^{*}_{r}(\mc G)$ is stably finite,
\item There are no Tarski matrices of any degree over $\mc G^{a}$,
\item Every clopen subset of $\mc G^{(0)}$ is not paradoxical with respect to $\mc G^{a}$
\item There exists a faithful invariant mean on $\mc G^{a}$ normalised at $1$,
\item There is a faithful $\mc G$-invariant measure on $\mc G^{(0)}$.
\end{enumerate}
\end{Thm}
\begin{proof}
We prove (i)$\implies$ (ii) $\implies$ (iii) $\implies$ (v) $\implies$ (i), since (iii) $\iff$ (iv) easily follows from \cite[Theorem 3.13]{MR3448402} and \cite[Theorem 5.11]{Bonicke-Li} and (v) $\iff$ (vi) is a consequence of \cite[Corollary 4.12]{MR3548134}.

(i)$\implies$(ii) is classical.
 %, and the proof can be found for instance in \cite{}.
(ii)$\implies$(iii) is the previous Corollary.
(iii)$\implies$(v) is the Tarski alternative of Kudryatseva--Lawson--Lenz--Resende \cite[Theorem 3.13]{MR3448402}.
(v)$\implies$(i) follows directly from \cite[Proposition 4.7]{MR3548134}.%, which just constructs the faithful trace directly. 
\end{proof}

We will extend this minimal result in the final sections of this text to introduce condition on the full group under some additional hypotheses on the groupoid.

\section{Polycyclic monoids, Cuntz algebras and Thompson's group $V$}\label{Sect:polycyclic-thompson}
We outline one important example of the general machinery we have developed in the previous sections; %as a prelude to adding critera to the classification results from the previous section
This is putting the cart before the horse, as historically this example is motivating for, and could be considered the inception of the theory of Boolean inverse monoids \cite{MR2317635}, \cite{MR2372319}  and their tight groupoids \cite{MR2419901}.

\begin{Def}
The polycyclic monoids $P_{n}$ are defined by the inverse monoid presentation $\langle a_{1},...,a_{n} \mid a_{i}^{*}a_{j}=\delta_{i,j}\rangle$, where $\delta_{i,j}$ is either a global identity if $i=j$ or a zero element if $i\neq j$. 

These inverse monoids arise naturally as a model of ``paradoxical behaviour" in geometry: suppose we have a representation of $P_{n}$ in the symmetric inverse monoid $I(X)$ for some infinite set $X$.  the condition ``$a_{i}^{*}a_{j}=\delta_{i,j}$" means that the domains of the $a_{i}$ are $X$, and the ranges are disjoint subsets of $X$.

A representation of $P_{n}$ is \textit{tight} %in the sense of Example \ref{Ex:semigroupgroupoid}
if the union of the ranges of the $a_{i}$ is the whole of the set $X$. A natural and concrete representation of $P_{n}$ is obtained on the binary tree.
\end{Def}

Polycyclic monoids are not themselves Boolean inverse monoids, but Lawson in \cite{MR2372319} associated a Boolean inverse monoid $C(P_{n})$ to $P_{n}$, which is universal for tight representations of $P_{n}$. This inverse monoid can be alternatively obtained as the Boolean inverse monoid of bisections of the transformation groupoid $\mc G_{tight}(P_{n})$, which is topologically principal and minimal. It is known that the full group $F(\mc G_{tight}(P_{n}))$ is isomorphic to Thompson's group $V_{n,1}$ \cite[Theorem 4.2]{Lawson-minimal}; in particular, for $n=2$ we recover Thompson's group $V$.%, and an illustration of how the isomorphism can be constructed is shown in Figure \ref{Fig:PolytoThompsons}\footnote{Not drawn yet}.

On the \Cst-algebraic side case, the tightness condition corresponds to the \textit{Cuntz relation} -- the sum of the range projections $a_{i}a_{i}^{*}$ is equal to the unit, thus the tight groupoid \Cst-algebra is isomorphic to the Cuntz algebra $\mc O_{n}$ \cite{MR2419901}. 

The representation $\pi: F(\mc G_{tight}(P_{n})) \rightarrow C_{c}(\mc G_{tight}(P_{n}))$ from subsection \ref{Sect:groupoidrep} was described independently by Nekrasheyvch \cite{MR2119267} and Birget \cite{MR2104771}, and was used explicitly in Haagerup--Olesen \cite{Haagerup-Olesen} to intertwine the $C^{*}$-simplicity of $F$ and the amenability of $F$.

\subsection{What do we know about the \Cst-algebras associated to $C(P_{2})$}

By inducing the representations of $V$ to $C(P_{2})$ we can prove the following.

\begin{Thm}
$C^{*}C(P_{2})$ is not exact, and $C^{*}_{r}(C(P_{2}))$ is not nuclear. 
\end{Thm}
\begin{proof}
Inducing the maximal representation $\pi$ of $V$ to $C(P_{2})$ we obtain a representation $\sigma:=\Ind_{V}^{C(P_{2})}\pi$ of $C(P_{2})$, the restriction $\Res_{V}^{C(P_{2})}(\sigma)$ contains $\pi$, and so is weakly equivalent to $\pi$ by its universality. Thus $C^{*}F_{2} \subset C^{*}V \subset C^{*}_{\sigma}(C(P_{2}))$, as $F_{2} \leqslant V$, hence $C^{*}_{\sigma}(C(P_{2}))$ is not exact because $C^{*}F_{2}$ is not exact and exactness passes to subalgebras. The first half of the result follows, as $C^{*}C(P_{2})$ quotients onto $C^{*}_{\sigma}(C(P_{2}))$, and exactness passes to quotients.

To see the second claim, recall that $C^{*}_{r}(C(P_{2}))$ is isomorphic to a groupoid \Cst-algebra, coming from the transformation groupoid associated with the natural action of $C(P_{2})$ on the filter space $\widehat{E}$. If $C^{*}_{r}(C(P_{2}))$ were nuclear, then this transformation groupoid would be measurewise amenable by \cite[Corollary 6.2.14]{MR1799683}.%\todo{Add an argument that $C(P_2)$ is a semilattice and use [Exel, Inverse semigroups and combinatorial C*-algebras])}.
However, by by \cite[Proposition 6.2.8]{MR1799683} this would imply that $C(P_{2})$ had the weak containment property, and so $C^{*}C(P_{2})$ would be nuclear. However, this is not the case, as it fails to be exact.
\end{proof}

We remark this leaves the possibility that $C^{*}_{r}C(P_{2})$ is exact open, and this would imply that $V$ was exact using Proposition \ref{prop:exact-subgroup-via-nuclearity}, since exactness of $C_r^*(V)$ follows from permanence properties \cite[Proposition 10.2.3, Theorem 10.2.5]{MR2391387} and characterisation of exact groups via exactness of the reduced $C^*$-algebra due to Kirchberg and Wassermann \cite[Theorem 5.1.10]{MR2391387}.

As was pointed out to us by the anonymous referee, this seems to be the first explicit example of an inverse semigroup with a faithful tight representation whose tight C*-algebra is nuclear but whose universal one is not. 

\section{Full group representations of Thompson's group $V$}

Let $X$ be a compact totally disconnected space and let $\mc G$ be a Hausdorff \'etale groupoid with base space $X$. We remark that we are not asking for $\mc G$ to be second countable here. The main result of this section is a combination of the results from \cite{MR3448402},  \cite{MR3548134} and \cite{Bonicke-Li} which can be seen as a unified Tarski alternative for ample groupoids, and we refer the reader to these sources for more information. We continue to use some definitions from the previous section.

Let us recall two of \textit{Green's relations} in the context of inverse semigroups:

\begin{Def}
Let $S$ be an inverse semigroup and let $e,f \in E(S)$ be its idempotents. Then:
\begin{itemize}
\item $e \mc D f$ if there is an $s \in S$ such that $e=s^{*}s$ and $f=ss^{*}$,
\item $e \leq_{\mc J} f$ if there exists $e'\in E(S), e' \leq f$ such that $e\mc D e'$ and
\item $e \mc J f$ if $e \leq_{\mc J} f$ and $f \leq_{\mc J} e$. 
\end{itemize}
\end{Def}
We remark that if $S=I(X)$ is the inverse semigroup of partial bijections on a set $X$ and $A,B\subseteq X$, then $1_{A}\leqslant_{\mc J} 1_{B}$ if there is an injection from $A$ to $B$, and $1_{A} \mc D 1_{B}$ if there is a bijection between $A$ and $B$. 

\begin{Thm}\label{Thm:generalgroupoid}
Let $\mc G$ be a topologically amenable, ample, principal Hausdorff groupoid with compact unit space $X$ such that $\mc G^{a}$ satisfies $\mc D=\mc J$ and has an almost unperforated type semigroup. Then the following are equivalent:
\begin{enumerate}
\item There is a Tarski matrix of degree 2 over $\mc G^{a}$.
\item There is a clopen subset $E \subset X$ such that $E$ is $(2,1)$-paradoxical.
\item $V \leqslant F(\mc G)$.
\end{enumerate}
\end{Thm}
\begin{proof}
(i)$\iff$(ii) follows from \cite[Theorem 3.13]{MR3448402} and \cite[Theorem 5.11]{Bonicke-Li}.

(i)$\implies$ (iii): This is essentially \cite[Remark 2.4]{MR3448402}. A Tarski matrix of degree $2$ over $\mc G^{a}$ gives us a weakly paradoxical clopen subset. But since $\mc D=\mc J$, we know that weak paradoxicality is the same as strong paradoxicality. From the discussion in Section \ref{Sect:polycyclic-thompson}, we know that a strongly paradoxical set induces a tight representation of $P_{2}$ in $\mc G^{a}$, thus a representation of $F(G_{tight}(P_{2}))=V$ into $F(\mc G)$.

(iii)$\implies$ (ii) is more complicated. Under the assumptions on $\mc G$, $(k,l)$-paradoxicality implies $(2,1)$-paradoxicality by the argument of \cite{MR3448402}. So we know by the Tarski alternatives that there is an $\mc G^{a}$-invariant mean supported on $E$, and thus there is a $\mc G$-invariant measure supported on $E$ for every clopen subset $E\subset X$.

Suppose that $V \leqslant F(\mc G)$. Then $V$ acts on $\mc G^{(0)}=X$ by homemorphisms through $F(\mc G)$ via the groupoid representation $\pi_{\kappa}$. Using Zorn's lemma, we extract a $V$-invariant minimal closed subsystem $Z$. As $Z$ is an intersection of clopen subsets $Z_{n}$ we can construct an invariant measure $\mu$ on $Z$: each $Z_{n}$ admits a $\mc G$-invariant probability measure $\mu_{n}$ using the discussion above and then we can take $\mu$ to be the weak* limit of these. It is $\mc G$-invariant, thus it is also $F(\mc G)$-invariant. Hence, using either the classification results of \cite{MR3231220} or \cite{LB-MB}, we can conclude that the action is topologically free.

The goal now is to show firstly that this action is actually free, as opposed to just topologically free, and then to embed the transformation groupoid it generates into a topologically amenable groupoid as a closed subgroupoid. We approach this in a general Lemma below.

\begin{Lemma}\label{Lemma:midproof}
Let $\mc G$ be a Hausdorff, ample groupoid such that $\Germs(\mc G)$ is a Hausdorff, principal groupoid and let $\mc H$ be a topologically principal, Hausdorff closed subgroupoid of $\mc G$. Then $\mc H$ is principal, and $\mc H$ is homeomorphic to a closed subgroupoid of $\Germs(\mc G)$.
\end{Lemma}
\begin{proof}
Suppose first that $\mc H$ is not principal, but only essentially principal. Then there exists an element $\gamma \in \mc H$ such that $\gamma \in \mc H_{s(\gamma)}^{s(\gamma)}$, and $\gamma \not = s(\gamma)$. 

As elements of $\mc G$ however $\gamma$ and $s(\gamma)$ must map to the same germ, as $\Germs(\mc G)$ is principal -- that is, for every clopen bisection $A\subset\mc G$ containing $\gamma$ satisfies $[A_{s(\gamma)}]=[s(A)_{s(\gamma)}]$. However, this implies that they must map to equivalent germs over $\mc H$, since $\mc H$ uses the subspace topology from $\mc G$. This is impossible as $\mc H$ is already an effective groupoid by Proposition \ref{prop:nek}.

Note that this also shows that the surjective map $\mc G \rightarrow \Germs(\mc G)$ from $\mc G$ to $\Germs(\mc G)$ is injective on $\mc H$, and it maps the bisections $A \subset \mc H$ onto themselves, as they are equal to their own germs. This proves that the map is continuous -- as the preimage of a set $O_{A}$ is necessarily just $A$, and both of these are clopen, and thus gives the desired homeomorphism.
\end{proof}

We apply this to the subgroupoid $\mc H = V \curvearrowright Z \subset F(\mc G) \curvearrowright X$. We observe $\mc H$ is closed, as the topology is inherited from the product topology, and topologically principal and so we can apply Lemma \ref{Lemma:midproof}, and the action of $V$ is therefore free, and $V \curvearrowright Z$ is embedded in $\Germs(F(\mc G) \curvearrowright X)$ as a closed subgroupoid. 

To complete the proof, observe that $\Germs(F(\mc G) \curvearrowright X)$ is a closed subgroupoid of $\mc G \cong \Germs(G^{a}\curvearrowright X)$ by Lemma \ref{Lem:germing} and noncommutative Stone duality \cite[Theorem 2.22]{Lawson-minimal}. As topological amenability passes to closed subgroupoids, $\Germs(F(\mc G) \curvearrowright X)$ is also topologically amenable. The same closed subgroupoid argument then implies $V \curvearrowright Z$ is topologically amenable. But the presence of an $V$-invariant probability measure on $Z$ implies that in fact $V$ is amenable. This is a contradiction, since $V$ has free subgroups.
\end{proof}

\begin{Cor}
Let $\mc G$ be an ample, minimal, topologically amenable Hausdorff groupoid with compact unit space, let $\mc G^{a}$ be the Boolean inverse monoid of bisections and let $F(\mc G)$ be the topological full group of $\mc G$. If $T(\mc G^{a})$ has no perforation,then the following are equivalent:
\begin{enumerate}
\item $C^{*}(\mc G)$ admits a faithful tracial state,
\item $C^{*}(\mc G)$ is stably finite,
\item There are no Tarski matrices of any degree over $\mc G^{a}$,
\item Every clopen subset of $\mc G^{(0)}$ is not paradoxical with respect to $\mc G^{a}$
\item There exists a faithful invariant mean on $\mc G^{a}$ normalised at $1$,
\item There is a faithful $\mc G$-invariant measure on $\mc G^{(0)}$.
\item There is no embedding of Thompson's group $V$ into $F(\mc G)$.
\item There is no embedding of $C(P_{2})$ into $\mc G^{a}$.
\end{enumerate}
\end{Cor}

\begin{Rem}
The above results are also connected to recent work by Matte Bon \cite{MB} considering general questions on embeddings of full groups.
\end{Rem}

\subsection{Applying Theorem \ref{Thm:generalgroupoid} to wobbling groups}

We begin this section with a few cursory definitions. Let $X$ be a uniformly discrete metric space of bounded geometry.

\begin{Def}
A partial bijection $t: A \rightarrow B$, $A,B \subseteq X$ is a \textit{partial translation} if there exists a $C>0$ such that for every $x \in A$ the distance $d(x,t(x))\leqslant C$. The collection of all partial translations, $I_{b}(X)$ is a Boolean inverse monoid with idempotent subalgebra $E(I_{b}(X)) = \textbf{2}^{X}$, its group of units $W(X)$ is called the wobbling group of $X$.
\end{Def}

\begin{Def} (The coarse groupoid)
Let $X$ be uniformly discrete with bounded geometry. Then $I_{b}(X)$ acts on $\beta X$, the Stone-\v{C}ech compactification of $X$, by partial homeomorphisms. The coarse groupoid, denoted by $G(X)$, is the transformation groupoid of this action.
\end{Def}

This definition is one of the original definitions of the coarse groupoid given by Skandalis-Tu-Yu in \cite{MR1905840}, which is the most convenient for the purposes of this paper. We collect a few properties of this in the next proposition.

\begin{Prop}\label{Prop:wobblinggroupisfull}
The wobbling group $W(X)$ is the full group of the ample, principal, Hausdorff topological groupoid $G(X)$.
\end{Prop}
\begin{proof}
The assumptions of this statement are collecting what is known from the original work of Skandalis--Tu--Yu \cite{MR1905840} concerning $G(X)$. There, $G(X)$ was constructed as the transformation groupoid of $I_{b}(X)$. It suffices to show that $G(X)^{a}=I_{b}(X)$, which follows by observing that any bisection $A$ induces a translation $t_{A}:s(A)\cap X \rightarrow r(A)\cap X$, the graph of which closes back to $A$.
\end{proof}

Note that applying the construction of the algebraic groupoid representation from Section \ref{Sect:groupoidrep} to the wobbling group, we obtain a representation of $W(X)$ in the unitaries of the uniform Roe algebra  $C^{*}_{r}(G(X))$. We know that this is exact as a \Cst-algebra whenever $X$ has property A \cite{MR1905840}, and admits a trace whenever $X$ is amenable. The following strengthening of property A implies the existence of a faithful trace. 

\begin{Def}
A uniformly discrete metric space of bounded geometry $X$ is \textit{supramenable} if no subset of $X$ is paradoxical.
\end{Def}

The following result is thus Theorem \ref{Thm:generalgroupoid} applied to coarse groupoids.

\begin{Thm}\label{Thm:wobbling}
Let $X$ be a uniformly discrete metric space of bounded geometry with property A, then the following are equivalent:
\begin{enumerate}
\item $X$ is not supramenable,
\item $V \leqslant W(X)$.
\end{enumerate}
\end{Thm}

We remark that there are representations of $V$ into wobbling groups of amenable groups, notably those of solvable groups with exponential growth (containing binary trees) as in the approach of \cite{MR3351042} for $\SL_{3}(\mathbb{Z})$.

\end{document}